\documentclass[12pt,a4paper]{article}
\usepackage{amsmath}
\usepackage{amssymb}
\usepackage{t1enc}
\usepackage[latin1]{inputenc}
\usepackage[english]{babel}
\usepackage{amsfonts}
\usepackage{latexsym}
\usepackage{bm}
\usepackage{setspace}
\usepackage[toc,page]{appendix}
\usepackage{graphicx}
\usepackage{enumerate}
\usepackage[numbers,sort&compress]{natbib}
\usepackage{tikz}

\usepackage{theoremref}
\usepackage{lineno}

\usepackage{mathtools}

\DeclarePairedDelimiter\floor{\lfloor}{\rfloor}

\newtheorem{theorem}{Theorem}
\newtheorem{prop}{Proposition}
\newtheorem{lemma}{Lemma}

\newtheorem{construction}{Construction}

\allowdisplaybreaks

\begin{document}
\title{Isolation of regular graphs and $k$-chromatic graphs
}

\author{Peter Borg \\[2mm]
\normalsize Department of Mathematics \\
\normalsize Faculty of Science \\
\normalsize University of Malta\\
\normalsize Malta\\
\normalsize \texttt{peter.borg@um.edu.mt}
}

\date{}
\maketitle

\begin{abstract}
Given a set $\mathcal{F}$ of graphs, we call a copy of a graph in $\mathcal{F}$ an $\mathcal{F}$-graph. The $\mathcal{F}$-isolation number of a graph $G$, denoted by $\iota(G,\mathcal{F})$, is the size of a smallest set $D$ of vertices of $G$ such that the closed neighbourhood of $D$ intersects the vertex sets of the $\mathcal{F}$-graphs contained by $G$ (equivalently, $G - N[D]$ contains no $\mathcal{F}$-graph). Thus, $\iota(G,\{K_1\})$ is the domination number of $G$. For any integer $k \geq 1$, let $\mathcal{F}_{1,k}$ be the set of regular graphs of degree at least $k-1$, let $\mathcal{F}_{2,k}$ be the set of graphs whose chromatic number is at least $k$, and let $\mathcal{F}_{3,k}$ be the union of $\mathcal{F}_{1,k}$ and $\mathcal{F}_{2,k}$. Thus, $k$-cliques are members of both $\mathcal{F}_{1,k}$ and $\mathcal{F}_{2,k}$. We prove that for each $i \in \{1, 2, 3\}$, $\frac{m+1}{{k \choose 2} + 2}$ is a best possible upper bound on $\iota(G, \mathcal{F}_{i,k})$ for connected $m$-edge graphs $G$ that are not $k$-cliques. The bound is attained by infinitely many (non-isomorphic) graphs. The proof of the bound depends on determining the graphs attaining the bound. This appears to be a new feature in the literature on isolation. Among the result's consequences are a sharp bound of Fenech, Kaemawichanurat and the present author on the $k$-clique isolation number and a sharp bound on the cycle isolation number. 
\\

\end{abstract}

\section{Introduction} \label{Introsection}
 
For standard terminology in graph theory, we refer the reader to \cite{West}. Most of the terminology used in this paper is defined in \cite{Borg}. 

The set of positive integers is denoted by $\mathbb{N}$. For $n \in \{0\} \cup \mathbb{N}$, $[n]$ denotes the set $\{i \in \mathbb{N} \colon i \leq n\}$. Note that $[0]$ is the empty set $\emptyset$. Arbitrary sets and graphs are taken to be finite. For a non-empty set $X$, ${X \choose 2}$ denotes the set of $2$-element subsets of $X$.

Every graph $G$ is taken to be \emph{simple}, meaning that its vertex set $V(G)$ and edge set $E(G)$ satisfy $E(G) \subseteq {V(G) \choose 2}$. An edge $\{v,w\}$ may be represented by $vw$. If $n = |V(G)|$, then $G$ is called an \emph{$n$-vertex graph}. If $m = |E(G)|$, then $G$ is called an \emph{$m$-edge graph}. For $v \in V(G)$, $N_{G}(v)$ denotes the set of neighbours of $v$ in $G$, $N_{G}[v]$ denotes the closed neighbourhood $N_{G}(v) \cup \{ v \}$ of $v$, and $d_{G}(v)$ denotes the degree $|N_{G} (v)|$ of $v$. For $X \subseteq V(G)$, $N_G[X]$ denotes the closed neighbourhood $\bigcup_{v \in X} N_G[v]$ of $X$, $G[X]$ denotes the subgraph of $G$ induced by $X$, and $G-X$ denotes the graph obtained by deleting the vertices in $X$ from $G$. Thus, $G[X] = (X,E(G) \cap {X \choose 2})$ and $G - X = G[V(G) \backslash X]$. With a slight abuse of notation, for $S \subseteq E(G)$, $G - S$ denotes the graph obtained by removing the edges in $S$ from $G$, that is, $G - S = (V(G), E(G) \backslash S)$. For $x \in V(G) \cup E(G)$, $G - \{x\}$ may be abbreviated to $G - x$. The subscript $G$ may be omitted where no confusion arises; for example, $N_G(v)$ may be abbreviated to $N(v)$. 

For $n \geq 1$, $K_n$ and $P_n$ denote the graphs $([n], {[n] \choose 2})$ and $([n], \{\{i,i+1\} \colon i \in [n-1]\})$, respectively. For $n \geq 3$, $C_n$ denotes $([n], \{\{1,2\}, \{2,3\}, \dots, \{n-1,n\}, \{n,1\}\})$. A copy of $K_n$ is called a \emph{complete graph} or an \emph{$n$-clique}. A copy of $P_n$ is called an \emph{$n$-path} or simply a \emph{path}. A copy of $C_n$ is called an \emph{$n$-cycle} or simply a \emph{cycle}. A $3$-cycle is a $3$-clique and is also called a \emph{triangle}. If $G$ contains a $k$-clique $H$, then we call $H$ a \emph{$k$-clique of $G$}.

If $D \subseteq V(G) = N[D]$, then $D$ is called a \emph{dominating set of $G$}. The \emph{domination number of $G$}, denoted by $\gamma(G)$, is the size of a smallest dominating set of $G$. If $\mathcal{F}$ is a set of graphs and $F$ is a copy of a graph in $\mathcal{F}$, then we call $F$ an \emph{$\mathcal{F}$-graph}. A subset $D$ of $V(G)$ is called an \emph{$\mathcal{F}$-isolating set of $G$} if no subgraph of $G-N[D]$ is an $\mathcal{F}$-graph. The \emph{$\mathcal{F}$-isolation number of $G$}, denoted by $\iota(G, \mathcal{F})$, is the size of a smallest $\mathcal{F}$-isolating set of $G$. We may abbreviate $\iota(G, \{F\})$ to $\iota(G, F)$. Clearly, $D$ is a dominating set of $G$ if and only if $D$ is a $\{K_1\}$-isolating set of $G$. Thus, $\gamma(G) = \iota(G, K_1)$.

Caro and Hansberg introduced the isolation problem in \cite{CaHa17}. It is a natural generalization of the classical domination problem \cite{C, CH, HHS, HHS2, HL, HL2}. Ore \cite{Ore} proved that $\gamma(G) \leq n/2$ for any connected $n$-vertex graph $G$ with $n \geq 2$ (see \cite{HHS}). This is one of the earliest results in domination theory. While the deletion of the closed neighbourhood of a dominating set produces the \emph{null graph} (the graph with no vertices), the deletion of the closed neighbourhood of a $\{K_2\}$-isolating set produces an \empty{empty graph} (a graph with no edges). Caro and Hansberg~\cite{CaHa17} proved that if $G$ is a connected $n$-vertex graph, then $\iota(G, K_2) \leq n/3$ unless $G$ is a $2$-clique or a $5$-cycle. Solving a problem of Caro and Hansberg~\cite{CaHa17}, Fenech, Kaemawichanurat and the present author~\cite{BFK} proved that 
\begin{equation} \iota(G, K_k) \leq \frac{n}{k+1} \label{BFKbound}
\end{equation} 
unless $G$ is a $k$-clique or $k = 2$ and $G$ is a $5$-cycle, and that the bound is sharp. Ore's result and the Caro--Hansberg result are the cases $k = 1$ and $k = 2$, respectively. Let $\mathcal{C}$ be the set of cycles. Solving another problem of Caro and Hansberg~\cite{CaHa17}, the present author~\cite{Borg} proved that 
\begin{equation} \iota(G,\mathcal{C}) \leq \frac{n}{4} \label{Borgbound}
\end{equation} 
unless $G$ is a triangle, and that the bound is sharp. Domination and isolation have been particularly investigated for maximal outerplanar graphs \cite{BK, BK2, CaWa13, CaHa17, Ch75, DoHaJo16, DoHaJo17, HeKa18, LeZuZy17, Li16, MaTa96, KaJi, To13}, mostly due to connections with Chv\'{a}tal's Art Gallery Theorem \cite{Ch75}.

Fenech, Kaemawichanurat and the present author~\cite{BFK2} also obtained a sharp upper bound on $\iota(G,K_k)$ in terms of the number of edges. The result is the analogue of (\ref{BFKbound}) given by Theorem~\ref{result1}. We need the following construction from \cite{BFK2}.  

\begin{construction}[\cite{BFK2}] \label{const2} \emph{Consider any $m, k \in \mathbb{N}$. Let $t_{k} = {k \choose 2} + 2$ and $q = \floor{ \frac{m+1}{t_k} }$. Then, $m+1 = qt_k + r$ for some $r \in \{0\} \cup [t_k - 1]$. Let $Q_{m,k}$ be a $q$-element set. If $q \geq 1$, then let $v_1, \dots, v_q$ be the elements of $Q_{m,k}$, let $U_1, \dots, U_q$ be $k$-element sets such that $U_1, \dots, U_q, Q_{m,k}$ are pairwise disjoint, and for each $i \in [q]$, let $u_i^1, \dots, u_i^k$ be the elements of $U_i$, and let $G_i = (U_i \cup \{v_i\}, {U_i \choose 2} \cup \{u_i^1v_i\})$. If either $q = 0$, $T = (\emptyset, \emptyset)$, and $G$ is an $m$-edge tree $T'$, or $q \geq 1$, $T$ is a tree with $V(T) = Q_{m,k}$, $T'$ is an $r$-edge tree with $V(T') \cap \bigcup_{i=1}^q V(G_i) = \{v_q\}$, and $G = \left( V(T') \cup \bigcup_{i=1}^q V(G_i), E(T) \cup E(T') \cup \bigcup_{i=1}^q E(G_i) \right)$, then we say that $G$ is an \emph{$(m,k)$-special graph} with \emph{quotient graph $T$} and \emph{remainder graph $T'$}, and for each $i \in [q]$, we call $G_i$ and $v_i$ a \emph{$k$-clique constituent of $G$} and the \emph{$k$-clique connection of $G_i$ in $G$}, respectively. We say that an $(m,k)$-special graph is \emph{pure} if $r = 0$. (Figure~1 in \cite{BFK2} is an illustration of a pure $(m, k)$-special graph.)}
\end{construction}

\begin{theorem}[\cite{BFK2}] \label{result1}
If $k \geq 1$ and $G$ is a connected $m$-edge graph that is not a $k$-clique, then
\begin{equation} \iota(G,K_k) \leq \frac{m+1}{{k\choose2}+2}. \label{cliquebnd}
\end{equation}
Moreover:\\
(i) Equality in (\ref{cliquebnd}) holds if and only if either $G$ is a pure $(m,k)$-special graph or $k=2$ and $G$ is a $5$-cycle. \\
(ii) If $G$ is an $(m,k)$-special graph, then $\iota(G,K_k) = \left \lfloor (m+1)/\left( {k\choose2}+2 \right ) \right \rfloor.$
\end{theorem}

Problem~7.4 in \cite{CaHa17} asks for bounds on $\iota(G,\mathcal{F})$ for other interesting sets $\mathcal{F}$. Section~1 of that seminal paper makes particular mention of $k$-colourable graphs. We address 
the problem for such graphs and the problem for regular graphs together. The main result presented here generalizes Theorem~\ref{result1} in these two desirable directions.

If $V(G) \neq \emptyset$ and $d(v) = r$ for each $v \in V(G)$, then $G$ is said to be \emph{$r$-regular} or simply \emph{regular}, and $r$ is called the \emph{degree of $G$}. For $k \geq 1$, let $\mathcal{F}_{1,k}$ be the set of regular graphs whose degree is at least $k-1$.

If there exists a function $f \colon V(G) \rightarrow [k]$ such that $f(v) \neq f(w)$ for every $v, w \in V(G)$ with $vw \in E(G)$, then $G$ is said to be \emph{$k$-colourable}, and $f$ is called a \emph{proper $k$-colouring of $G$}. The \emph{chromatic number of $G$}, denoted by $\chi(G)$, is the smallest non-negative integer $k$ such that $G$ is $k$-colourable. If $k = \chi(G)$, then $G$ is said to be \emph{$k$-chromatic}. For $k \geq 1$, let $\mathcal{F}_{2,k}$ be the set of graphs whose chromatic number is at least $k$. 
 
Let $\mathcal{F}_{3,k}$ be the union of $\mathcal{F}_{1,k}$ and $\mathcal{F}_{2,k}$. In Section~\ref{Proofsection1}, we prove the following result.

\begin{theorem} \label{mainresult2} If $\ell \in \{1, 2, 3\}$, $k \geq 1$, and $G$ is a connected $m$-edge graph that is not a $k$-clique, then 
\begin{equation} \iota(G, \mathcal{F}_{\ell,k}) \leq \frac{m+1}{{k\choose2}+2}. \label{thmbnd}
\end{equation}
Moreover:\\
(i) Equality in (\ref{thmbnd}) holds if and only if $G$ is a pure $(m,k)$-special graph, or $k=2$ and $G$ is a $5$-cycle, or $\ell \in \{1,3\}$, $k = 3$, and $G$ is a $4$-cycle. \\
(ii) If $G$ is an $(m,k)$-special graph, then $\iota(G, \mathcal{F}_{\ell,k}) = \left \lfloor (m+1)/\left( {k\choose2}+2 \right ) \right \rfloor.$
\end{theorem}
The proof of the bound depends on determining the graphs attaining the bound (i.e. the proof of (\ref{thmbnd}) uses (i) and (ii)). This appears to be a new feature in the literature on isolation.

Since $k$-cliques are $(k-1)$-regular and cycles are $2$-regular, Theorem~\ref{result1} and the following analogue of (\ref{Borgbound}) are immediate consequences.

\begin{theorem} \label{corcycle}
If $G$ is a connected $m$-edge graph that is not a triangle, then
\begin{equation} \iota(G,\mathcal{C}) \leq \frac{m+1}{5}. \label{cyclebnd}
\end{equation}
Moreover:\\
(i) Equality in (\ref{cyclebnd}) holds if and only if $G$ is a pure $(m,3)$-special graph or a $4$-cycle. \\
(ii) If $G$ is an $(m,3)$-special graph, then $\iota(G,\mathcal{C}) = \left \lfloor (m+1)/5 \right \rfloor.$
\end{theorem}

\section{Proof of Theorem~\ref{mainresult2}} \label{Proofsection1}

In this section, we prove Theorem~\ref{mainresult2}. We start with two lemmas from \cite{Borg}.

\begin{lemma}[\cite{Borg}] \label{lemma}
If $G$ is a graph, $\mathcal{F}$ is a set of graphs, $X \subseteq V(G)$, and $Y \subseteq N[X]$, then \[\iota(G, \mathcal{F}) \leq |X| + \iota(G-Y, \mathcal{F}).\] 
\end{lemma}

The set of components of a graph $G$ will be denoted by ${\rm C}(G)$.

\begin{lemma}[\cite{Borg}] \label{lemmacomp} 
If $G$ is a graph and $\mathcal{F}$ is a set of graphs, then \[\iota(G,\mathcal{F}) = \sum_{H \in {\rm C(G)}} \iota(H,\mathcal{F}).\]
\end{lemma}

A basic result in graph theory is that for any connected graph $G$, $|E(G)| \geq |V(G)|-1$, and equality holds if and only if $G$ is a tree.

\begin{prop} \label{edgespecialgraph}
If $G$ is a pure $(m,k)$-special graph with exactly $q$ $k$-clique constituents, then $m = ({k\choose2}+2)q-1$ and $\iota(G,\mathcal{F}_{\ell,k}) = q$ for each $\ell \in \{1, 2, 3\}$.
\end{prop}
\textbf{Proof.} Suppose that $G$ is a pure $(m,k)$-special graph with exactly $q$ $k$-clique constituents as in Construction~\ref{const2}. Then, $\{v_1, \dots, v_q\}$ is an $\mathcal{F}_{\ell,k}$-isolating set of $G$, so $\iota(G,\mathcal{F}_{\ell,k}) \leq q$. The subgraphs $G[U_1], \dots, G[U_q]$ of $G$ are $k$-cliques, which are $(k-1)$-regular and $k$-chromatic. Thus, if $D$ is an $\mathcal{F}_{\ell,k}$-isolating set of $G$, then $D \cap V(G_i) \neq \emptyset$ for each $i \in [q]$. Therefore, $\iota(G,\mathcal{F}_{\ell,k}) = q$. Now clearly $m = ({k\choose2}+1)q + |E(T)|$. Since $T$ is a $q$-vertex tree, $|E(T)| = q-1$. Thus, $m = ({k\choose2}+2)q-1$. \hfill{$\Box$}
\\

If $G$ is a copy of a graph $H$, then we write $G \simeq H$. For a graph $G$, the maximum degree of $G$ (that is, $\max\{d(v) \colon v \in V(G)\}$) is denoted by $\Delta(G)$.

We will use the following classical result.

\begin{theorem}[Brooks' Theorem \cite{Brooks}] \label{Brooksthm} If $G$ is a connected $n$-vertex graph, then $\chi(G) \leq \Delta(G)$ unless $G \simeq K_n$ or $n$ is odd and $G \simeq C_n$.
\end{theorem}

\noindent
\textbf{Proof of Theorem~\ref{mainresult2}.}  Let $t_k = {k \choose 2} + 2$ and $\mathcal{F} = \mathcal{F}_{3,k}$.

The argument in the proof of Proposition~\ref{edgespecialgraph} yields (ii). If $G$ is a $4$-cycle, then $\iota(G,\mathcal{F}_{\ell,1}) = 2 < \frac{5}{2} = \frac{m+1}{t_1}$, $\iota(G,\mathcal{F}_{\ell,2}) = 1 < \frac{m+1}{t_2}$, $\iota(G,\mathcal{F}_{1,3}) = \iota(G,\mathcal{F}_{3,3}) = 1 = \frac{m+1}{t_3}$, $\iota(G,\mathcal{F}_{2,3}) = 0 < \frac{m+1}{t_3}$, and if $k \geq 4$, then $\iota(G,\mathcal{F}_{\ell,k}) = 0 < \frac{m+1}{t_k}$. If $G$ is a $5$-cycle, then $\iota(G,\mathcal{F}_{\ell,1}) = 2 < \frac{6}{2} = \frac{m+1}{t_1}$, $\iota(G,\mathcal{F}_{\ell,2}) = 2 = \frac{m+1}{t_2}$, $\iota(G,\mathcal{F}_{\ell,3}) = 1 < \frac{m+1}{t_3}$, and if $k \geq 4$, then $\iota(G,\mathcal{F}_{\ell,k}) = 0 < \frac{m+1}{t_k}$. Together with Proposition~\ref{edgespecialgraph}, this settles the sufficiency condition in (i). 

Using induction on $m$, we now prove that the bound in (\ref{thmbnd}) holds, and we prove that it is attained only if $G$ and $k$ are as in (i). For the bound, since $\iota(G, \mathcal{F})$ is an integer and $\iota(G,\mathcal{F}_{\ell,k}) \leq \iota(G,\mathcal{F})$, it suffices to prove that $\iota(G, \mathcal{F}) \leq \frac{m+1}{t_k}$. If $k \leq 2$, then a $\{K_k\}$-isolating set of $G$ is an $\mathcal{F}_{\ell,k}$-isolating set of $G$, so the result is given by Theorem~\ref{result1}. Consider $k \geq 3$. The result is trivial if $m \leq 2$ or $\iota(G,\mathcal{F}) = 0$. Suppose that $m \geq 3$, $\iota(G,\mathcal{F}) \geq 1$, and $G$ is not a $k$-clique. Let $n = |V(G)|$. Since $G$ is connected, $m \geq n - 1$, so $n \leq m + 1$.\\

\noindent
\textbf{Case 1}: \emph{$G$ has no $k$-cliques.} \\

\noindent
\textbf{Subcase 1.1}: \emph{$d(v) \leq k-2$ for some $v \in V(G)$.} Let $G' = G-v$. Let $H_1, \dots, H_r$ be the distinct components of $G'$. 

Consider any $i \in [r]$. Let $D_i$ be a smallest $\mathcal{F}$-isolating set of $H_i$. Let $H_i' = H_i - N[D_i]$, let $c_i = \chi(H_i')$, and let $f_i$ be a proper $c_i$-colouring of $H_i'$. We have that $c_i \leq k-1$ and that $H_i'$ contains no regular graph of degree at least $k-1$. Since $G$ has no $k$-cliques, $|D_i| \leq \frac{|E(H_i)|+1}{t_k}$ by the induction hypothesis.  

Let $D = \bigcup_{i=1}^r D_i$. Suppose that $G - N[D]$ contains a connected regular graph $R$ of degree at least $k-1$. Since $d_R(v) \leq d_G(v)\leq k-2$, $v \notin R$. We obtain that for some $i \in [r]$, $H_i'$ contains $R$, a contradiction. Thus, $G - N[D]$ contains no regular graph of degree at least $k-1$. Let $f' \colon V(G'-N_{G'}[D]) \rightarrow [k-1]$ such that $f'(w) = f_i(w)$ for each $i \in [r]$ and each $w \in V(H_i')$. Since $d(v) \leq k-2$, $[k-1]$ has an element $j$ that is not in $\{f'(w) \colon w \in N(v) \cap V(G' - N_{G'}[D])\}$. Let $f \colon V(G-N[D]) \rightarrow [k-1]$ such that $f(w) = f'(w)$ for each $w \in V(G'-N_{G'}[D])$, and such that if $v \in V(G-N[D])$, then $f(v) = j$. Since $f$ is a proper $(k-1)$-colouring of $G - N[D]$, $D$ is an $\mathcal{F}$-isolating set of~$G$.

Since $G$ is connected, for each $i \in [r]$, $vw_i \in E(G)$ for some $w_i \in V(H_i)$. We have 
\[ m \geq |\bigcup_{i=1}^r (\{vw_i\} \cup E(H_i))| \geq \sum_{i=1}^r t_k|D_i| = t_k|D| \geq t_k \iota(G,\mathcal{F}),\]
so $\iota(G,\mathcal{F}) < \frac{m+1}{t_k}$. \\

\noindent
\textbf{Subcase 1.2}: \emph{$d(v) \geq k-1$ for each $v \in V(G)$.} Then, $\Delta(G) \geq k-1$. Let $v \in V(G)$ with $d(v) = \Delta(G)$. Let $G' = G - N[v]$. 

Suppose $d(v) = k-1$. Then, $G$ is $(k-1)$-regular. Since $G$ is not a $k$-clique, we have $n \geq k+1$, so $V(G') \neq \emptyset$. Suppose that $G'$ contains a $(k-1)$-regular graph $R'$. Since $G$ is connected, $uw \in E(G)$ for some $u \in V(R')$ and some $w \in V(G) \backslash V(R')$. We have $d_G(u) \geq d_{R'}(u) + 1 \geq k$, which contradicts $\Delta(G) = k-1$. Thus, $G'$ contains no $(k-1)$-regular graph. Since $\Delta(G') \leq \Delta(G) = k-1$, $\chi(G') \leq k-1$ by Theorem~\ref{Brooksthm}. Thus, $\{v\}$ is an $\mathcal{F}$-isolating set of $G$. By the handshaking lemma, we have $2m = (k-1)n \geq (k-1)(k+1) = 2t_k + k - 5$, so $m + 1 \geq t_k + \frac{k-3}{2} \geq t_k$. We have $\iota(G, \mathcal{F}) = 1 \leq \frac{m+1}{t_k}$. Suppose $\iota(G, \mathcal{F}) = \frac{m+1}{t_k}$. Then, $n = k+1$ and $k = 3$. Thus, $G$ is a $4$-cycle.

Now suppose $d(v) \geq k$. Let $s = d(v)$. Suppose $N[v] = V(G)$. By the handshaking lemma, 
\[2m = d(v) + \sum_{x \in N(v)} d(x) \geq d(v) + d(v)(k-1) \geq k + k(k-1) = 2t_k + k-4 \geq 2t_k - 1, \]
so $t_k < m+1$. Since $\{v\}$ is an $\mathcal{F}$-isolating set of $G$, $\iota(G,\mathcal{F}) = 1 < \frac{m+1}{t_k}$. Now suppose $N[v] \neq V(G)$. Then, $V(G') \neq \emptyset$. Let $H_0 = G[N(v)]$, and let $H_1, \dots, H_r$ be the distinct components of $G'$. Since $G$ is connected, for each $i \in [r]$, $x_iy_i \in E(G)$ for some $x_i \in N(v)$ and $y_i \in V(H_i)$. Let $A = \{xy \in E(G) \colon x \in N(v), \, y \in V(G')\}$. Then, $x_1y_1, \dots, x_ry_r \in A$, so $|A| \geq r$. By the handshaking lemma, 
\begin{align} 2|E(H_0)| &= \sum_{x \in N(v)} d_{H_0}(x) = \sum_{x \in N(v)} |N(x) \backslash (\{v\} \cup V(G'))| \nonumber \\
&= \sum_{x \in N(v)} (d(x) - 1 - |N(x) \cap V(G')|) \nonumber \\
&\geq (k-2)d(v) - \sum_{x \in N(v)} |N(x) \cap V(G')| = (k-2)s - |A|, \nonumber
\end{align}
so $|E(H_0)| \geq \frac{(k-2)s - |A|}{2}$. Let $A' = \{vx \colon x \in N(v)\} \cup E(H_0) \cup A$. Then, 
\begin{align} |A'| &= s + |E(H_0)| + |A| \geq \frac{sk + |A|}{2} = \frac{(k-1)k}{2} + \frac{(s-k+1)k + |A|}{2} \nonumber \\
&= t_k + \frac{(s-k+1)k + |A| - 4}{2}. \nonumber
\end{align}
Since $E(G) = A' \cup \bigcup_{i=1}^r E(H_i)$ and $|A| \geq r$, 
\begin{equation} m = |A'| + \sum_{i=1}^r |E(H_i)| \geq t_k + \frac{(s-k+1)k + r - 4}{2} + \sum_{i=1}^r |E(H_i)|. \label{m_ineq1.2}
\end{equation}
For each $i \in [r]$, let $D_i$ be a smallest $\mathcal{F}$-isolating set of $H_i$. Let $D = \{v\} \cup \bigcup_{i=1}^r D_i$. Then, $D$ is an $\mathcal{F}$-isolating set of $G$. Since $G$ has no $k$-cliques, for each $i \in [r]$, $H_i$ is not a pure $(|E(H_i)|,k)$-special graph. 

Suppose $k \neq 3$ or $H_i \not\simeq C_4$ for each $i \in [r]$. By the induction hypothesis, for each $i \in [r]$, $|D_i| < \frac{|E(H_i)|+1}{t_k}$, and hence, since $t_k|D_i| < |E(H_i)| + 1$, $t_k|D_i| \leq |E(H_i)|$. By~(\ref{m_ineq1.2}), 
\[m \geq t_k \left( 1 + \sum_{i=1}^r |D_i| \right) + \frac{(s-k+1)k + r - 4}{2} \geq t_k|D| + \frac{k + r - 4}{2} \geq t_k |D|\] 
as $k \geq 3$ and $r \geq 1$. We have $\iota(G, \mathcal{F}) \leq |D| < \frac{m+1}{t_k}$.

Now suppose $k = 3$ and $H_j \simeq C_4$ for some $j \in [r]$.  We may assume that $j = 1$. Let $z_1$ and $z_2$ be the two members of $N_{H_1}(y_1)$, and let $z_3$ be the member of $V(H_1) \backslash N_{H_1}[y_1]$. We have $V(H_1) = \{y_1, z_1, z_2, z_3\}$. Let $G_1 = G - \{y_1, z_1, z_2\}$. Note that $N_{G_1}(v) = N_G(v)$. If $d_{G_1}(z_3) = 0$, then the components of $G_1$ are $(\{z_3\}, \emptyset)$ and $G_1 - z_3$, and we take $I$ to be $G_1 - z_3$. If $d_{G_1}(z_3) \geq 1$, then $G_1$ is connected, and we take $I$ to be $G_1$. Since $d_{I}(v) = d_G(v) \geq k = 3$, $I$ is not a $4$-cycle. Since $G$ has no $k$-cliques, $I$ is not a pure $(|E(I)|,k)$-special graph. By the induction hypothesis, $\iota(I, \mathcal{F}) < \frac{|E(I)|+1}{t_k} = \frac{|E(I)|+1}{5}$. Let $D_I$ be an $\mathcal{F}$-isolating set of $I$ of size $\iota(I, \mathcal{F})$. Since $\{y_1, z_1, z_2\} \subset N[y_1]$, $\{y_1\} \cup D_I$ is an $\mathcal{F}$-isolating set of $G$. Let $J = \{x_1y_1, y_1z_1, y_1z_2, z_1z_3, z_2z_3\}$. Note that $J \subseteq E(G) \backslash E(I)$. We have
\[\iota(G,\mathcal{F}) \leq 1 + |D_I| < 1 + \frac{|E(I)|+1}{5} = \frac{|E(I)| + |J| + 1}{5} \leq \frac{m+1}{5} = \frac{m+1}{t_k}. \]
\vskip 2mm

\noindent
\textbf{Case 2}: \emph{$G$ has at least one $k$-clique.} Let $F$ be a $k$-clique of $G$. Since $G$ is connected and $G \not\simeq K_k$, $N[v] \backslash V(F) \neq \emptyset$ for some $v \in V(F)$. Let $u \in N[v] \backslash V(F)$. We have $V(F) \subseteq N[v]$, so $|N[v]| \geq k+1$ and
\begin{equation}\label{eq 0}
|E(G[N[v]])| \geq t_k - 1.
\end{equation}

Suppose $V(G) = N[v]$. Then, $\{v\}$ is an $\mathcal{F}$-isolating set of $G$, so $\iota(G,\mathcal{F}) = 1 \leq \frac{m+1}{t_k}$ by (\ref{eq 0}). If $\iota(G,\mathcal{F}) = \frac{m+1}{t_k}$, then $V(G) = V(F) \cup \{u\}$ and $E(G) = E(F) \cup \{uv\}$, so $G$ is a pure $(m,k)$-special graph.

Now suppose $V(G) \neq N[v]$. Let $G' = G-N[v]$ and $n' = |V(G')|$. We have
\begin{equation}\label{eq 2}
n \geq n' + k + 1
\end{equation}
and $n' \geq 1$. Let $\mathcal{H} = {\rm C}(G')$. For any $H \in \mathcal{H}$ and $x \in N(v)$ such that $xy \in E(G)$ for some $y \in V(H)$, we say that $H$ is \emph{linked to $x$}. Since $G$ is connected, for each $H \in \mathcal{H}$, $H$ is linked to at least one member of $N(v)$, so $x_Hy_H \in E(G)$ for some $x_H \in N(v)$ and some $y_H \in V(H)$. Clearly,
\begin{equation}\label{edgesboundineqsharp1}
m \geq |E(F) \cup \{ uv \}| + \sum_{H \in \mathcal{H}} |E(H) \cup \{x_Hy_H\}| = t_k - 1 + \sum_{H \in \mathcal{H}} (|E(H)| + 1).
\end{equation}

Let
\[\mathcal{H}' = \{H \in \mathcal{H} \colon H \simeq K_k\}.\]
By the induction hypothesis, for each $H \in \mathcal{H} \backslash \mathcal{H}'$, $\iota(H,\mathcal{F}) \leq \frac{|E(H)| + 1}{t_k}$, and equality holds only if $H$ is a pure $(|E(H)|,k)$-special graph or $k = 3$ and $H$ is a $4$-cycle.
\\

\noindent \textbf{Subcase 2.1:} $\mathcal{H}' = \emptyset$. By Lemma~\ref{lemma} (with $X = \{v\}$ and $Y = N[v]$), Lemma~\ref{lemmacomp}, and (\ref{edgesboundineqsharp1}),
\begin{align}\label{edgesboundineqsharp3}
\iota(G,\mathcal{F}) &\leq 1 + \iota(G',\mathcal{F}) = 1 + \sum_{H \in \mathcal{H}} \iota(H,\mathcal{F}) \leq 1 + \sum_{H \in \mathcal{H}} \frac{|E(H)|+1}{t_k} \leq \frac{m+1}{t_k}.
\end{align}
Suppose $\iota(G,\mathcal{F}) = \frac{m+1}{t_k}$. Then, the inequality in (\ref{edgesboundineqsharp1}) is an equality, and $\iota(H,\mathcal{F}) = \frac{|E(H)| + 1}{t_k}$ for each $H \in \mathcal{H}$. 

Let $H_1, \dots, H_r$ be the distinct members of $\mathcal{H}$. Let $I = \{i \in [r] \colon H_i$ is a pure $(|E(H_i)|,k)$-special graph$\}$. For each $i \in I$, let $G_{i,1}, \dots, G_{i,q_i}$ be the distinct $k$-clique constituents of $H_i$, and for each $j \in [q_i]$, let $v_{i,j}$ be the $k$-clique connection of $G_{i,j}$ in $H_i$, and let $u_{i,j}^1, \dots, u_{i,j}^k$ be the vertices of $G_{i,j} - v_{i,j}$, where $u_{i,j}^1$ is the neighbour of $v_{i,j}$ in $G_{i,j}$. By Proposition~\ref{edgespecialgraph}, $|E(H_i)| = t_kq_i - 1$ for each $i \in I$. If $i \in [r] \backslash I$, then, by the above, $k = 3$, $H_i$ is a $4$-cycle, and $|E(H_i)| + 1 = t_k$. Since equality holds throughout in (\ref{edgesboundineqsharp1}), we have 
\begin{equation} V(G) = V(F) \cup \{u\} \cup \bigcup_{i=1}^r V(H_i) \label{compare2.1-1}
\end{equation}
and
\begin{equation} E(G) = E(F) \cup \{uv\} \cup \{x_{H_1}y_{H_1}, \dots, x_{H_r}y_{H_r}\} \cup \bigcup_{i=1}^r E(H_i). \label{compare2.1-2}
\end{equation}
Let 
\[D = \{v\} \cup \{y_{H_i} \colon i \in [r] \backslash I\} \cup \left( \bigcup_{i \in I} \bigcup_{j=1}^{q_i} \{u_{i,j}^1\} \right).\]
Then, $D$ is an $\mathcal{F}$-isolating set of $G$, and 
\[|D| = 1 + \left(\sum_{i \in [r] \backslash I} \frac{|E(H_i)| + 1}{t_k} \right)  + \sum_{i \in I} q_i = 1 + \sum_{i=1}^{r} \frac{|E(H_i)|+1}{t_k} = \frac{m+1}{t_k} = \iota(G,\mathcal{F}).\]

Suppose that $[r] \backslash I$ has an element $j$. Then, $k = 3$ and $H_j$ is a $4$-cycle. Suppose $x_{H_j} = u$. Then, we have that $(D \backslash \{v, y_{H_j}\}) \cup \{u\}$ is an $\mathcal{F}$-isolating set of $G$ of size $|D|-1$, contradicting $|D| = \iota(G,\mathcal{F})$. Thus, $x_{H_j} \in V(F)$. This yields that $D \backslash \{v\}$ is an $\mathcal{F}$-isolating set of $G$, contradicting $|D| = \iota(G,\mathcal{F})$. 

Therefore, $[r] \backslash I = \emptyset$, and hence $I = [r]$. Consider any $i \in [r]$. Suppose $x_{H_i} \in V(F)$ and $y_{H_i} = v_{i,j}$ for some $j \in [q_i]$. Then, the components of $G - N[(D \backslash \{v, u_{i,j}^1\}) \cup \{v_{i,j}\}]$ are $G_{i,j} - \{v_{i,j}, u_{i,j}^1\}$ and $G[\{u\} \cup V(F-x_{H_i})]$, and clearly they contain no $\mathcal{F}$-graphs (a proper $(k-1)$-colouring of $G[\{u\} \cup V(F-x_{H_i})]$ is obtained by assigning different colours to the $k-1$ vertices of $F - x_{H_i}$, and colouring $u$ differently from $v$). Thus, we have that $(D \backslash \{v, u_{i,j}^1\}) \cup \{v_{i,j}\}$ is an $\mathcal{F}$-isolating set of $G$ of size $|D| - 1$, a contradiction. Suppose $x_{H_i} \in V(F)$ and $y_{H_i} = u_{i,j}^l$ for some $j \in [q_i]$ and $l \in [k]$. As above, we obtain that $(D \backslash \{v, u_{i,j}^1\}) \cup \{u_{i,j}^l\}$ is an $\mathcal{F}$-isolating set of $G$ of size $|D| - 1$, a contradiction. Thus, $x_{H_i} = u$. Suppose $y_{H_i}= u_{i,j}^l$ for some $j \in [q_i]$ and $l \in [k]$. Then, we have that $(D \backslash \{v,u_{i,j}^1\}) \cup \{u\}$ is an $\mathcal{F}$-isolating set of $G$ of size $|D| - 1$, a contradiction. Thus, $y_{H_i}= v_{i,j}$ for some $j \in [q_i]$. Since this has been established for each $i \in [r]$, we have that $G[\{u\} \cup \bigcup_{i=1}^r \bigcup_{j=1}^{q_i} \{v_{i,j}\}]$ is a tree, and hence $G$ is a pure $(m,k)$-special graph with $G[N[v]], G_{1,1}, \dots, G_{1,q_1}, \dots, G_{r,1}, \dots, G_{r,q_r}$ being its $k$-clique constituents and $u$ being the $k$-clique connection of $G[N[v]]$ in $G$. 
\\

\noindent \textbf{Subcase 2.2:} $\mathcal{H}' \neq \emptyset$. Let $H' \in \mathcal{H}'$. Then, $H'$ is linked to some $x \in N(v)$. Let
\vskip 8 pt

\indent $\mathcal{H}_{1} = \{H \in \mathcal{H}' \colon H$ is only linked to $x\}$ and
\vskip 8 pt

\indent $\mathcal{H}_{2} = \{H \in \mathcal{H} \backslash \mathcal{H}' \colon H$ is only linked to $x\}$.
\vskip 8 pt

\noindent 
For each $H \in \mathcal{H}_2$, let $D_H$ be a smallest $\mathcal{F}$-isolating set of $H$. Then, $|D_H| = \iota(H,\mathcal{F}) \leq \frac{|E(H)|+1}{t_k}$ by the induction hypothesis.
\\

\noindent
\textbf{Subcase 2.2.1:} \emph{Each member of $\mathcal{H}'$ is linked to at least two members of $N(v)$.} Thus, $\mathcal{H}_1 = \emptyset$. Let 
\[\mbox{$X = \{x\} \cup V(H')$ \quad and \quad $G^* = G - X$.}\] 
Then, a component $G^*_v$ of $G^*$ satisfies $N[v] \backslash \{x\} \subseteq V(G^*_v)$, and the members of $\mathcal{H}_2$ are the other components of $G^*$. Let $D^*_v$ be a smallest $\mathcal{F}$-isolating set of $G^*_v$. Since $H'$ is linked to $x$ and to some $x' \in N(v) \backslash \{x\}$, there exist $y, y' \in V(H')$ with $xy, x'y' \in E(G)$. 
Let $D = D^*_v \cup \{y\} \cup \bigcup_{H \in \mathcal{H}_2} D_H$. Since $X \subseteq N[y]$, $D$ is an $\mathcal{F}$-isolating set of $G$. Thus,
\begin{equation}
\iota(G,\mathcal{F}) \leq |D^*_v| + 1 + \sum_{H \in \mathcal{H}_2} |D_H| \leq |D^*_v| + \frac{t_k}{t_k} + \sum_{H \in \mathcal{H}_2} \frac{|E(H)| + 1}{t_k}. \nonumber 
\end{equation}
We have
\begin{align}
m &\geq |E(G^*_v) \cup \{vx\}| + |E(H') \cup \{xy, x'y'\}| + \sum_{H \in \mathcal{H}_2} |E(H) \cup \{x_Hy_H\}| \nonumber \\
  &= |E(G^*_v)| + 1 + t_k + \sum_{H \in \mathcal{H}_2} (|E(H)| +1). \label{edgesineq1}
\end{align}

If $G^*_v \not\simeq K_k$, then $|D^*_v| \leq \frac{|E(G^*_v)|+1}{t_k}$ by the induction hypothesis, so 
\[\iota(G,\mathcal{F}) \leq \frac{|E(G^*_v)| + 1}{t_k} + \frac{t_k}{t_k} + \sum_{H \in \mathcal{H}_2} \frac{|E(H)| + 1}{t_k} < \frac{m+1}{t_k}. \]

Now suppose $G^*_v \simeq K_k$. Since $|N[v] \backslash \{x\}| \geq k$ and $N[v] \backslash \{x\} \subseteq V(G^*_v)$, $V(G^*_v) = N[v] \backslash \{x\}$. Let $Y = (X \cup V(G^*_v)) \backslash \{v,x,y\}$ and $G_Y = G - \{v,x,y\}$. Then, the components of $G[Y]$ and the members of $\mathcal{H}_2$ are the components of $G_Y$.

Suppose that $G[Y]$ contains no $\mathcal{F}$-graph. Since $v, y \in N[x]$, $\{x\} \cup \bigcup_{H \in \mathcal{H}_2} D_H$ is an $\mathcal{F}$-isolating set of $G$, so $\iota(G,\mathcal{F}) \leq 1 + \sum_{H \in \mathcal{H}_2} |D_H|$. Since $G^*_v \simeq K_k$, (\ref{edgesineq1}) gives us $m + 1 \geq 2t_k + \sum_{H \in \mathcal{H}_2} (|E(H)|+1) > t_k(1 + \sum_{H \in \mathcal{H}_2} |D_H|) \geq t_k\iota(G,\mathcal{F})$, so $\iota(G,\mathcal{F}) < \frac{m+1}{t_k}$. 

Now suppose that $G[Y]$ contains an $\mathcal{F}$-graph $F_Y$. Then, $\chi(F_Y) \geq k$ or $F_Y$ is a regular graph of degree at least $k-1$. Thus, $|N_{G[Y]}[z]| \geq k$ for some $z \in V(F_Y)$; this is given by Theorem~\ref{Brooksthm} if $\chi(F_Y) \geq k$. Let $W \subseteq N_{G[Y]}[z]$ such that $z \in W$ and $|W| = k$. Let $G_1 = G_v^*$, $G_2 = H'$, $v_1 = v$, $v_2 = y$, $G_1' = G_1 - v_1$, and $G_2' = G_2 - v_2$. We have
\begin{equation}\label{edgesmain.2}
N_{G[Y]}[z] \subseteq Y = V(G_1') \cup V(G_2').
\end{equation}
Since $|V(G_1')| = |V(G_2')| = k-1$, it follows that $|W \cap V(G_1')| \geq 1$ and $|W \cap V(G_2')| \geq 1$. By (\ref{edgesmain.2}), $z \in V(G_j')$ for some $j \in \{1, 2\}$. Let $Z = V(G_j) \cup W$. Since $z \in V(G_j)$ and $G_j \simeq K_k$,
\begin{equation} Z \subseteq N[z]. \label{edgesmain.3}
\end{equation}
We have
\begin{equation} |Z| = |V(G_j)| + |W \backslash V(G_j)| = k + |W \cap V(G_{3-j}')| \geq k+1. \label{edgesmain.4}
\end{equation}
Let $G_Z = G - Z$. Then, $V(G_Z) = \{x\} \cup (V(G_{3-j}) \backslash W) \cup \bigcup_{H \in \mathcal{H}_2} V(H)$. The components of $G_Z - x$ are the clique $G_Z[V(G_{3-j}) \backslash W]$ (which has less than $k$ vertices) and the members of $\mathcal{H}_2$. Moreover, $v_{3-j} \in V(G_{3-j}) \backslash W$ (by (\ref{edgesmain.2})), $v_{3-j} \in N_{G_Z}(x)$, and $N_{G_Z}(x) \cap V(H) \neq \emptyset$ for each $H \in \mathcal{H}_2$ (by the definition of $\mathcal{H}_2$). Thus, $G_Z$ is connected.

Suppose $\mathcal{H}_2 \neq \emptyset$. Then, $G_Z \not\simeq K_k$. 
By (\ref{edgesmain.3}), Lemma~\ref{lemma}, and the induction hypothesis, $\iota(G,\mathcal{F}) \leq 1 + \iota(G_Z,\mathcal{F}) \leq 1 + \frac{|E(G_Z)|+1}{t_k}$. Let $w^* \in W \cap V(G_{3-j}')$ (recall that $|W \cap V(G_{3-j}')| \geq 1$). Since $G_{3-j} \simeq K_k$, $w^*w \in E(G_{3-j})$ for each $w \in V(G_{3-j}) \backslash \{w^*\}$. We have 
\begin{align} m  &\geq |E(G_j)| + |\{xv_j, zw^*\}| + |\{w^*w \colon w \in V(G_{3-j}) \backslash \{w^*\}\}| + |E(G_Z)| \nonumber \\
&\geq {k \choose 2} + 1 + k + |E(G_Z)| \geq t_k + 2 + |E(G_Z)| \nonumber
\end{align}
as $k \geq 3$. Thus, since $\iota(G,\mathcal{F}) \leq 1 + \frac{|E(G_Z)|+1}{t_k}$, $\iota(G,\mathcal{F}) < \frac{m+1}{t_k}$.

Now suppose $\mathcal{H}_2 = \emptyset$. Then, $G^* = G^*_v$, so $V(G) = V(G^*_v) \cup \{x\} \cup V(H')$. Recall that the $k$-clique $H'$ is linked to at least two members of $N(v)$. Thus, 
\[\mbox{$n = 2k+1$ \quad and \quad $m \geq 2{k \choose 2} + 3 = 2t_k - 1$.}\]
If $|N[w]| \geq k+2$ for some $w \in V(G)$, then $|V(G-N[w])| \leq k-1$, so $\{w\}$ is an $\mathcal{F}$-isolating set of $G$, and hence $\iota(G,\mathcal{F}) = 1 < \frac{m+1}{t_k}$. Suppose $|N[w]| \leq k+1$ for each $w \in V(G)$. Then, $\Delta(G) = d(v) = k$, $N[z] = Z = V(G_j) \cup \{w\}$ for some $w \in V(G_{3-j}')$ (by (\ref{edgesmain.3}) and (\ref{edgesmain.4})), and $V(G-N[z]) = \{x\} \cup V(G_{3-j} - w)$. If $G-N[z]$ contains no $\mathcal{F}$-graph, then $\iota(G,\mathcal{F}) = 1 < \frac{m+1}{t_k}$. Suppose that $G-N[z]$ contains an $\mathcal{F}$-graph $F'$. Since $|V(G-N[z])| = k$, $G-N[z] = F' \simeq K_k$. Since $\Delta(G) = k$, we have $N(x) = \{v_{j}\} \cup V(G_{3-j} - w)$ and, since $z \in V(G_j)$ and $w \in N[z] \cap V(G_{3-j})$, $N[w] = \{z\} \cup V(G_{3-j})$. Thus, $V(G-N[w]) = \{x\} \cup (V(G_j) \backslash \{z\})$. Since $|V(G-N[w])| = k \geq 3$ and $N[x] \cap V(G_j') = \emptyset$, $\{w\}$ is an $\mathcal{F}$-isolating set of $G$, so $\iota(G,\mathcal{F}) = 1 < \frac{m+1}{t_k}$.
\\

\noindent
\textbf{Subcase 2.2.2:} \emph{Some member of $\mathcal{H}'$ is linked to only one member of $N(v)$.} Thus, we may assume that $H'$ is linked to $x$ only. For each $H \in \mathcal{H}_1 \cup \mathcal{H}_2$, $y_H \in N(x)$ for some $y_H \in V(H)$. By the same argument in the corresponding Subcase~2.2.2 of the proof of Theorem~1.3 in \cite{BFK2}, we obtain that $\iota(G,\mathcal{F}) \leq \frac{m+1}{t_k}$, and that equality holds only if
\[E(G) = E(H') \cup \{xy_{H'}\} \cup \bigcup_{H \in \mathcal{H}_2'} (E(H) \cup \{xy_H\}),\]
where $\mathcal{H}_2' = \mathcal{H}_2 \cup \{G_v^*\}$ for a connected subgraph $G_v^*$ of $G$ with 
\[y_{G_v^*} = v \in V(G_v^*) = V(G) \backslash \left(\{x\} \cup V(H') \cup \bigcup_{H \in \mathcal{H}_2} V(H) \right),\] 
and $\iota(H,\mathcal{F}) = \frac{|E(H)| + 1}{t_k}$ for each $H \in \mathcal{H}_2'$. Thus, if $\iota(G,\mathcal{F}) = \frac{m+1}{t_k}$, then $G$ is as in the case $\iota(G,\mathcal{F}) = \frac{m+1}{t_k}$ in Subcase~2.1 (with $H'$ and $x$ now taking the role of $F$ and $u$, respectively, in (\ref{compare2.1-1}) and (\ref{compare2.1-2})), so $G$ is a pure $(m,k)$-special graph.~\hfill{$\Box$}
\\

\noindent
\textbf{Acknowledgements.} The author is grateful to the anonymous referee for checking the paper and providing constructive remarks.


\begin{thebibliography}{} 

\bibitem{Borg} P. Borg, Isolation of cycles, Graphs and Combinatorics 36 (2020), 631--637.

\bibitem{BFK} P. Borg, K. Fenech and P. Kaemawichanurat, Isolation of $k$-cliques, Discrete Mathematics 343 (2020), paper 111879.

\bibitem{BFK2} P. Borg, K. Fenech and P. Kaemawichanurat, Isolation of $k$-cliques II, Discrete Mathematics 345 (2022), paper 112641.

\bibitem{BK} P. Borg and P. Kaemawichanurat, Partial domination of maximal outerplanar graphs, Discrete Applied Mathematics 283 (2020), 306--314.

\bibitem{BK2} P. Borg and P. Kaemawichanurat, Extensions of the Art Gallery Theorem, Annals of Combinatorics 27 (2023) 31--50.

\bibitem{Brooks} R.L. Brooks, On colouring the nodes of a network, Proceedings of the Cambridge Philosophical Society 37 (1941), 194--197.

\bibitem{CaWa13} C.N. Campos and Y. Wakabayashi, On dominating sets of maximal outerplanar graphs, Discrete Applied Mathematics 161 (2013), 330--335.

\bibitem{CaHa17} Y. Caro and A. Hansberg, Partial domination - the isolation number of a graph, Filomat 31:12 (2017), 3925--3944.

\bibitem{Ch75} V. Chv\'{a}tal, A combinatorial theorem in plane geometry, Journal of Combinatorial Theory Series B 18 (1975), 39--41.

\bibitem{C} E.J. Cockayne, Domination of undirected graphs -- A survey, Lecture Notes in Mathematics, Volume 642, Springer, 1978, 141--147.

\bibitem{CH} E.J. Cockayne and S.T. Hedetniemi, Towards a theory of domination in graphs, Networks 7 (1977), 247--261.

\bibitem{DoHaJo16} M. Dorfling, J.H. Hattingh and E. Jonck, Total domination in maximal outerplanar graphs II, Discrete Applied Mathematics 339 (2016), 1180--1188.

\bibitem{DoHaJo17} M. Dorfling, J.H. Hattingh and E. Jonck, Total domination in maximal outerplanar graphs, Discrete Applied Mathematics 217 (2017), 506--511.

\bibitem{HHS} T.W. Haynes, S.T. Hedetniemi and P.J. Slater, Fundamentals of Domination in Graphs, Marcel Dekker, Inc., New York, 1998.

\bibitem{HHS2} T.W. Haynes, S.T. Hedetniemi and P.J. Slater (Editors), Domination in Graphs: Advanced Topics, Marcel Dekker, Inc., New York, 1998.

\bibitem{HL} S.T. Hedetniemi and R.C. Laskar (Editors), Topics on Domination, in: Annals of Discrete Mathematics, vol. 48, North-Holland Publishing Co., Amsterdam, 1991, Reprint of Discrete Mathematics 86 (1990), no. 1--3.

\bibitem{HL2} S.T. Hedetniemi and R.C. Laskar, Bibliography on domination in graphs and some basic definitions of domination parameters, Discrete Mathematics 86 (1990), 257--277.

\bibitem{HeKa18} M. A. Henning and P. Kaemawichanurat, Semipaired domination in maximal outerplanar graphs, Journal of Combinatorial Optimization 38 (2019), 911--926.

\bibitem{LeZuZy17} M. Lema\'{n}ska, R. Zuazua and P. \.{Z}yli\'{n}ski, Total dominating sets in maximal outerplanar graphs, Graphs and Combinatorics 33 (2017), 991--998.

\bibitem{Li16} Z. Li, E. Zhu, Z. Shao and J. Xu, On dominating sets of maximal outerplanar and planar graphs, Discrete Applied Mathematics 198 (2016), 164--169.

\bibitem{MaTa96} L. R. Matheson and R. E. Tarjan, Dominating sets in planar graphs, European Journal of Combinatorics 17 (1996), 565--568.

\bibitem{Ore} O. Ore, Theory of Graphs, American Mathematical Society Colloquium Publications, Volume 38, American Mathematical Society, Providence, R.I., 1962.

\bibitem{To13} S. Tokunaga, Dominating sets of maximal outerplanar graphs, Discrete Applied Mathematics 161 (2013), 3097--3099.

\bibitem{KaJi} S. Tokunaga, T. Jiarasuksakun and P. Kaemawichanurat, Isolation number of maximal outerplanar graphs, Discrete Applied Mathematics 267 (2019), 215--218.

\bibitem{West} D. B. West, Introduction to Graph Theory, second edition, Prentice Hall, 2001. 

\end{thebibliography}
\end{document}